\documentclass[12pt]{article}
\usepackage{mathrsfs}
\usepackage{amsmath,amssymb}
\newcommand{\be}{\begin{equation}}
\newcommand{\ee}{\end{equation}}
\newcommand{\bal}{\begin{aligned}}
\newcommand{\eal}{\end{aligned}}
\newcommand{\bee}{\begin{equation*}}
\newcommand{\eee}{\end{equation*}}
\usepackage{color}

\def\cwedge{\bigcirc\kern-1.07em\wedge\ }

\newtheorem{thm}{Theorem}[section]
\newtheorem{lem}[thm]{Lemma}

\newtheorem{question}{Question}

\numberwithin{equation}{section}
\begin{document}
\begin{center}
{\large \bf A remark concerning universal curvature
identities on 4-dimensional Riemannian manifolds}\\
\end{center}
\footnotetext{\small{\it E-mail addresses}: {\bf
prettyfish@skku.edu} (Y. Euh), {\bf chohee1108@skku.edu} (C. Jeong),
and {\bf parkj@skku.edu} (J. Park)}
\begin{center}
Yunhee Euh, Chohee Jeong, and JeongHyeong Park
\end{center}

{\small
\begin{center}
     Department of Mathematics,
    Sungkyunkwan University,
    Suwon 440-746, KOREA
\end{center}
}



\begin{abstract}
 We shall prove the universality of
the curvature identity for the 4-dimensional Riemannian manifold
using a different method than that used by Gilkey, Park, and
Sekigawa \cite{GPS}.
\end{abstract}


\section{Introduction}

Berger \cite{Beg} derived a curvature identity on a 4-dimensional
compact oriented Riemannian manifold $M=(M,g)$ from the generalized
Gauss-Bonnet formula
\begin{equation*}\label{eq:4GB}
 32\pi^2\chi(M)= \int_M \tau^2 -4 |\rho|^2 + |R|^2 dv,
  \end{equation*}
where $R$ is the curvature tensor, $\rho$ is the Ricci tensor and
$\tau$ is the scalar curvature of $M$.
 The curvature
identity is the quadratic equation which involves only the curvature
tensor and not its covariant derivatives as follows:
   \begin{equation}\label{main equation}
    \frac{1}{4}(|R|^2-4|\rho|^2+\tau^2)g-\check{R}+2\check{\rho}+L\rho-\tau\rho=0.
    \end{equation}
    Here,
     \begin{equation*}
     \begin{gathered}
     \check{R}:\check{R}_{ij}=\sum_{a,b,c}R_{abci}{R^{abc}}_{j},\qquad
     \check{\rho}:\check{\rho}_{ij}=\sum_a\rho_{ai}\rho^{a}_{~j},\\
     L:(L\rho)_{ij}=2\sum_{a,b}R_{iabj}\rho^{ab}.
     \end{gathered}
    \end{equation*}
\vskip0.3cm \noindent Euh, Park, and Sekigawa \cite{EPS} proved that
Equation \eqref{main equation} holds on the space of all Riemannian
metrics on any 4-dimensional Riemannian manifold, and gave some
applications of the curvature identity \cite{EPS2, EPS3}. Labbi
\cite{La} showed the same phenomena occurs for the higher
dimensional cases by using purely algebraic computations in the ring
of double forms and also provided some applications of the curvature
identity in \cite{La(2010)}. Recently, Gilkey, Park, and Sekigawa
\cite{GPS} gave a new proof of the curvature identity using heat
trace methods. Here, we raise the following question:

\begin{question}
Is there another curvature identity such as the quadratic curvature
identity \eqref{main equation} which holds on any 4-dimensional
Riemannian manifold $(M,g)$?
\end{question}

 In the present paper, we shall
give an answer to the above Question with a different method given
by \cite{GPS}. Namely, we shall prove the following
theorem.\vskip0.2cm

\begin{thm}
The curvature identity \eqref{main equation} is universal as a
symmetric 2-form valued quadratic curvature identity for a
4-dimensional Riemannian manifold.
\end{thm}

\vskip0.2cm {{The authors would like to express their thanks to
Professors P. Gilkey and K. Sekigawa for their helpful comments and
valuable suggestions.}}

\section{Preliminary}\label{sec3}

Let $M$ be an $m$-dimensional Riemannian manifold and
$\mathcal{I}_{m,n}^2$($n$ is even) be the space of symmetric
$2$-form valued invariants which are homogeneous of degree $n$ in
the derivatives of the metric on $M$. In \cite{GPS}, Gilkey, Park,
and Sekigawa proved that the universality of the curvature identity
in the setting of the space $\mathcal{I}_{4,4}^2$.
 Now, we set
  \begin{equation*}
    \begin{gathered}
    \Phi_1 : = |R|^2g,\ \ \ \Phi_2
:=|\rho|^2g,\ \ \  \Phi_3:=\tau^2g,\ \ \ \Phi_4 :=\check{R},\ \ \
\Phi_5
:=\check{\rho},\\
  \ \Phi_6:=L\rho, \ \ \Phi_7 :=\tau\rho, \ \  \Phi_8 = (\triangle\tau)g, \quad\Phi_9 = {\text{Hess }} \tau,\quad
\Phi_{10} = \tilde\triangle\rho,
 \end{gathered}
    \end{equation*}
where $\tilde\triangle\rho$ denotes the rough Laplacian acting on
the Ricci tensor $\rho$, namely locally expressed by
$(\tilde\triangle\rho)_{ij} = \sum_{a}\nabla^a \nabla_a \rho_{ij}$.
Then, we have the following:

\noindent
\begin{lem}\cite{GPS}
 \begin{enumerate}
\item $\mathcal{I}_{m,0}^2 = $ Span $\{g\}$,
\item $\mathcal{I}_{m,2}^2 = $ Span $\{\tau g, \rho \}$,
\item $\mathcal{I}_{m,4}^2 = $ Span $\{ \Phi_1, \Phi_2, \cdots, \Phi_7, \Phi_8, \Phi_9, \Phi_{10} \}$
\end{enumerate} \vskip0.5cm
\end{lem}


In \cite{GPS, GPS2}, Gilkey et al. proved that the curvature
identity

\begin{equation}\label{ci}
 \frac{\lambda}{4} \Phi_1-\lambda \Phi_2+\frac{\lambda}{4} \Phi_3-\lambda
 \Phi_4+2\lambda \Phi_5+\lambda \Phi_6-\lambda
 \Phi_7=0
\end{equation}
for any constant $\lambda(\neq0)$, is the only universal curvature
identity of this form if $m=4$ (\cite{GPS}, Theorem 1.2 (3) and
Lemma 1.4 (2)). We may easily check that the curvature identities
\eqref{main equation} and \eqref{ci} are equivalent to each other.
We emphasize that the invariance theory established by H. Weyl plays
an important role in their proof of the Theorem 1.2
\cite{GPS}.\\

Here, we give another direct proof for the same result by using
several test Riemannian manifolds of dimension 4.

\section{Proof of Main theorem}

We assume that the equality
 \begin{equation}\label{sum}
\sum_{i=1}^{10} c_i\Phi_i=0\
\end{equation}
holds for all 4-dimensional Riemannian manifolds. To prove Main
Theorem, it is sufficient to prove that $c_1=\frac{\lambda}{4}$,
$c_2=\lambda$, $c_3=\frac{\lambda}{4}$, $c_4=-\lambda$,
$c_5=2\lambda$, $c_6=\lambda$, $c_7=-\lambda$, $c_8=c_9=c_{10}=0$.

Applying \eqref{sum} to the test manifolds in Cases I, II, III, IV
and V, we will determine the coefficients $c_i$'s such that
$\sum_ic_i\Phi_i=0\ \ (i=1,\cdots, 10)$ by applying the method of
universal examples.
 This is the way we can show whether
the curvature identity \eqref{main equation} is universal or
not.\vskip0.2cm

\noindent{\bf Case I.}  Let $M$ be a locally product of Riemannian
surfaces $M^2(a)$ and $ M^2(b)$ of nonzero constant Gaussian
curvatures $a$ and $b$. Let $\{e_1, e_2\}$ and $\{e_3, e_4\}$ be the
orthonormal basis of $M^2(a)$ and $M^2(b)$, respectively. Then we
have the following:
\begin{equation}\label{eq:caseI}
\begin{gathered}
\Phi_1=4(a^2+b^2)I, \quad \Phi_2=2(a^2+b^2)I,\quad \Phi_3=4(a+b)^2
I,
 \\
\Phi_4=\begin{pmatrix}
  2a^2  & 0      & 0        & 0     \\
  0     & 2a^2   & 0        & 0     \\
  0     & 0      & 2b^2     & 0     \\
  0     & 0      & 0        & 2b^2
 \end{pmatrix}, \qquad
\Phi_5=\begin{pmatrix}
  a^2   & 0      & 0        & 0      \\
  0     & a^2    & 0        & 0\\
  0     & 0      & b^2      & 0 \\
  0     & 0      & 0        & b^2
   \end{pmatrix},\\
 \Phi_6=\begin{pmatrix}
  2a^2    & 0     & 0        & 0      \\
  0      & 2a^2   & 0        & 0      \\
  0      & 0     & 2b^2      & 0      \\
  0      & 0     & 0        & 2b^2
   \end{pmatrix}, \qquad
 \Phi_7=2(a+b)\begin{pmatrix}
  a     & 0      & 0        & 0      \\
  0     & a      & 0        & 0      \\
  0     & 0      & b        & 0      \\
  0     & 0      & 0        & b
 \end{pmatrix},\\
 {{\Phi_8 = \Phi_9 =\Phi_{10}=0.}}
\end{gathered}
\end{equation}
From \eqref{eq:caseI}, we can get two different equations such that
$\sum_{i}c_i\Phi_i=0$ :
\\
(I-i) (1,1)-component (or (2,2)-component)
\begin{equation*}
(4c_1+2c_2+4c_3+2c_4+c_5+2c_6+2c_7)a^2+(8c_3+2c_7)ab+(4c_1+2c_2+4c_3)b^2=0.
\end{equation*}
(I-ii) (3,3)-component (or (4,4)-component)
\begin{equation*}
(4c_1+2c_2+4c_3)a^2+(8c_3+2c_7)ab+(4c_1+2c_2+4c_3+2c_4+c_5+2c_6+2c_7)b^2=0.
\end{equation*}
We set $c_7=-\lambda$. Then from (I-i) and (I-ii), we have the
following relations:
\begin{equation}\label{eq:caseI_relation}
\begin{aligned}
& c_3=\frac{1}{4}\lambda,\\
& 4c_1+2c_2=-\lambda,\\
& 2c_4+c_5+2c_6=2\lambda.
\end{aligned}
\end{equation}

\noindent{\bf Case II.}  Let $M$ be a product of 3-dimensional
Riemannian manifold $M^3(a)$ of nonzero constant sectional curvature
$a$ and a real line $\mathbb R$. Let $\{e_1, e_2, e_3\}$ be the
orthonormal basis of $M^3(a)$. Then we have the following:
\begin{equation}\label{eq:caseII}
\begin{gathered}
\Phi_1=12a^2I,\qquad \Phi_2=12a^2I,\qquad \Phi_3=36a^2I,\\
\Phi_4=4a^2\begin{pmatrix}
  1      & 0      & 0       & 0     \\
  0      & 1      & 0       & 0     \\
  0      & 0      & 1       & 0     \\
  0      & 0      & 0       & 0
 \end{pmatrix},\quad
\Phi_5=4a^2\begin{pmatrix}
  1      & 0      & 0       & 0     \\
  0      & 1      & 0       & 0     \\
  0      & 0      & 1       & 0     \\
  0      & 0      & 0       & 0
 \end{pmatrix},\\
\Phi_6=8a^2\begin{pmatrix}
  1      & 0      & 0       & 0     \\
  0      & 1      & 0       & 0     \\
  0      & 0      & 1       & 0     \\
  0      & 0      & 0       & 0
 \end{pmatrix}, \quad
\Phi_7=12a^2\begin{pmatrix}
  1      & 0      & 0       & 0     \\
  0      & 1      & 0       & 0     \\
  0      & 0      & 1       & 0     \\
  0      & 0      & 0       & 0
 \end{pmatrix},
 \\
 {{\Phi_8 = \Phi_9 =\Phi_{10}=0.}}
\end{gathered}
\end{equation}
From \eqref{eq:caseII}, we can get two different equations such that
$\sum_{i}c_i\Phi_i=0$ :\\
(II-i) (1,1)-component ((2,2) or (3,3)-component)
\begin{equation*}
(3c_1+3c_2+9c_3+c_4+c_5+2c_6+3c_7)a^2=0.
\end{equation*}
(II-ii) (4,4)-component
\begin{equation*}
(c_1+c_2+3c_3)a^2=0.
\end{equation*}
From (II-i) and (II-ii), we have the following relation:
\begin{equation*}
\begin{aligned}
 c_4+c_5+2c_6+3c_7=0,
\end{aligned}
\end{equation*}
and hence, since $c_7=-\lambda$, we get
\begin{equation}\label{eq:caseII_relation}
\begin{aligned}
 c_4+c_5+2c_6=3\lambda.
\end{aligned}
\end{equation}
From \eqref{eq:caseI_relation} and \eqref{eq:caseII_relation}, we
have
    \begin{equation}\label{eq:re_caseII}
    \begin{aligned}
    c_4=-\lambda,\quad c_5+2c_6=4\lambda.
    \end{aligned}
    \end{equation}

\noindent{\bf Case III.} Let $M=M^4(a)$ be a space form of nonzero
constant sectional curvature $a$. Then we have the following:
\begin{equation}\label{eq:caseIII}
\begin{gathered}
 \Phi_1=24a^2I,\quad\Phi_2=36a^2 I,\quad
\Phi_3=144a^2I,\\
\Phi_4=6a^2I,\quad\ \Phi_5=9a^2I,\quad\ \Phi_6=18a^2I,\\
 {{\Phi_7=36a^2I,\ \quad\Phi_8 = \Phi_9 =\Phi_{10}=0.}}
\end{gathered}
\end{equation}
From \eqref{eq:caseIII}, we can get an equation such that
$\sum_{i}c_i\Phi_i=0$ :\\
(III) (1,1)-component ((2,2), (3,3), or (4,4)-component)
\begin{equation*}
(24c_1+36c_2+144c_3+6c_4+9c_5+18c_6+36c_7)a^2=0.
\end{equation*}
From (III-i), we have the following relation:
\begin{equation*}\label{eq:caseIII_relation}
8c_1+12c_2+48c_3+2c_4+3c_5+6c_6+12c_7=0.
\end{equation*}
Since $c_7=-\lambda$, from \eqref{eq:caseI_relation} and
\eqref{eq:re_caseII}, we get
    \begin{equation}\label{eq:re_caseIII}
    c_1=\frac{\lambda}{4}, \ \ c_2=-\lambda.
    \end{equation}

\noindent{\bf Case IV.} (\cite{EPS2}, Example 3.7) Let
$\mathfrak{g}=\text{span}_{\mathbb{R}}\{e_1,e_2,e_3,e_4\}$ be a
4-dimensional real Lie algebra equipped with the following Lie
bracket operation:
    \begin{equation}\label{eq:LieBracket}
    \begin{aligned}
    &[e_1,e_2]=ae_2,\qquad[e_1,e_3]=-ae_3-be_4,\qquad[e_1,e_4]=be_3-ae_4,\\
    &[e_2,e_3]=0,\quad\qquad[e_2,e_4]=0,\qquad\qquad\qquad[e_3,e_4]=0,
    \end{aligned}
    \end{equation}
where $a(\ne0)$, $b$ are constant. We define an inner product $<,>$
on $\mathfrak{g}$ by $<e_i,e_j>=\delta_{ij}$.  Let $G$ be a
connected and simply connected solvable Lie group with the Lie
algebra $\mathfrak{g}$ of $G$ and $g$ the $G$-invariant Riemannian
metric on $G$ determined by $<,>$. From \eqref{eq:LieBracket}, by
direct calculations, we have
    \begin{equation}
    \begin{gathered}
    R_{1212}=a^2,~~\quad R_{1313}=a^2,~~\quad R_{1414}=a^2,\\
    R_{2323}=-a^2,\quad R_{2424}=-a^2,\quad R_{3434}=a^2,
    \end{gathered}
    \end{equation}
and otherwise being zero up to sign.\\
\begin{equation*}
\begin{gathered} (\rho)=\begin{pmatrix}
  -3a^2     & 0     & 0     & 0      \\
  0     & a^2     & 0     & 0      \\
  0     & 0     & -a^2     & 0      \\
  0     & 0     & 0     & -a^2
 \end{pmatrix},   \quad \tau=-4a^2.
 \end{gathered}
\end{equation*}

\noindent Then, we have the following:
\begin{equation}\label{eq:caseIV}
\begin{gathered}
\Phi_1=24a^4I,\quad\Phi_2=12a^4I,\quad \Phi_3=16a^4I, \quad\Phi_4=6a^4I,\\
\Phi_5=a^4\begin{pmatrix}
  9     & 0     & 0     & 0      \\
  0     & 1     & 0     & 0      \\
  0     & 0     & 1     & 0      \\
  0     & 0     & 0     & 1
 \end{pmatrix},\quad
\Phi_6=2a^4\begin{pmatrix}
  1     & 0     & 0     & 0      \\
  0     & 1     & 0     & 0      \\
  0     & 0     & 5     & 0      \\
  0     & 0     & 0     & 5
 \end{pmatrix}, \\
 \Phi_7=4a^4\begin{pmatrix}
  3      & 0      & 0    & 0      \\
  0      & -1     & 0    & 0\\
  0      & 0      & 1    & 0 \\
  0      & 0      & 0    & 1
 \end{pmatrix},\quad
 \Phi_{10}=a^4\begin{pmatrix}
  8      & 0      & 0    & 0      \\
  0      & -8     & 0    & 0\\
  0      & 0      & -4    & 0 \\
  0      & 0      & 0    & -4
 \end{pmatrix}, \\
 \Phi_8 = \Phi_9 =0.
\end{gathered}
\end{equation}
{{From \eqref{eq:caseIV},  we can get three different equations such that $\sum_{i}c_i\Phi_i=0$ :\\
(IV-i) (1,1)-component\\
     \begin{equation}\label{1,1}
    (24c_1+12c_2+16c_3+6c_4+9c_5+2c_6+12c_7 {{+ 8c_{10}}})a^4=0.
    \end{equation}
(IV-ii) (2,2)-component
\begin{equation}\label{2,2}
    (24c_1+12c_2+16c_3+6c_4+c_5+2c_6-4c_7 {{- 8c_{10}}})a^4=0.
    \end{equation}
(IV-iii) (3,3)-component (or (4,4)-component)
\begin{equation}\label{3,3}
    (24c_1+12c_2+16c_3+6c_4+c_5+10c_6+4c_7 {{- 4c_{10}}})a^4=0.
\end{equation}
Thus, from \eqref{1,1}, taking account of \eqref{eq:caseI_relation},
\eqref{eq:re_caseII}, \eqref{eq:re_caseIII} and $a\ne 0$, we have
     \begin{equation}\label{560}
    -20\lambda +9c_5+2c_6 {{+ 8c_{10}}}=0.
    \end{equation}
Thus, from \eqref{2,2}, we have
 \begin{equation}\label{560-2}
    -4\lambda + c_5+2c_6+{{- 8c_{10}}}=0.
    \end{equation}
    Then, from \eqref{560} and \eqref{560-2}, we have
 \begin{equation}\label{56}
    5c_5+2c_6 = 12\lambda .
    \end{equation}
 Thus, from \eqref{eq:re_caseII} and \eqref{56},
  we have
\begin{equation}\label{562}
    c_5 = 2\lambda, \quad c_6 = \lambda.
    \end{equation}
    Thus, \eqref{560} and \eqref{562}, we have
\begin{equation}\label{10}
    c_{10} = 0.
    \end{equation}
\noindent{\bf Case V.} Let $M$ be the Riemannian product of
Riemannian surfaces $(M_1,g_1)$ and $(M_2,g_2)$, where the
Riemannian metrics
$g_1$ and $g_2$ are given locally by\\
\begin{equation*}
(g_1)=\begin{pmatrix}
  e^{2\sigma_1}     & 0          \\
  0     & e^{2\sigma_1}           \\
 \end{pmatrix}, \qquad \sigma_1=x_1^2+x_2^2
\end{equation*}
and
\begin{equation*}
(g_2)=\begin{pmatrix}
  e^{2\sigma_2}     & 0          \\
  0     & e^{2\sigma_2}           \\
 \end{pmatrix}, \qquad \sigma_2=x_3^2+x_4^2.
\end{equation*}\\
We set
\begin{equation*}
e_1=\frac{1}{e^{\sigma_1}}\frac{\partial}{\partial x_1}, \quad
e_2=\frac{1}{e^{\sigma_1}}\frac{\partial}{\partial x_2}, \quad
e_3=\frac{1}{e^{\sigma_2}}\frac{\partial}{\partial x_3}, \quad
e_4=\frac{1}{e^{\sigma_2}}\frac{\partial}{\partial x_4}.
\end{equation*}
We denote by $K_1$ and $K_2$ the Gaussian curvatures of $(M_1 , g_1
)$ and $(M_2 , g_2 )$, respectively. Then we have
\begin{equation} \label{gg}
K_1 =-4e^{-2\sigma_1},\quad K_2 =-4e^{-2\sigma_2}. \end{equation}
Thus, from \eqref{gg}, we have the scalar curvature
\begin{equation*}
\tau=-8e^{-2\sigma_1}-8e^{-2\sigma_2}.
\end{equation*}
Finally, we have\\
\begin{equation}\label{eq:caseV}
\begin{gathered}
\Phi_8=-64\big(e^{-4\sigma_1}(2\sigma_1-1)+e^{-4\sigma_2}(2\sigma_2-1)\big)I,\quad\Phi_9=\begin{pmatrix}
A   & 0     \\
0   & B
 \end{pmatrix},
\end{gathered}
\end{equation}
where
\begin{equation*}
\begin{gathered}
A=-32e^{-4\sigma_1}
    \begin{pmatrix}
    6x_1^2-2x_2^2-1  &  8x_1x_2\\
    8x_1x_2       &  -2x_1^2+6x_2^2-1
    \end{pmatrix},\\
    B=-32e^{-4\sigma_1}
    \begin{pmatrix}
    6x_3^2-2x_4^2-1  &  8x_3x_4\\
    8x_3x_4       &  -2x_3^2+6x_4^2-1
    \end{pmatrix}.
\end{gathered}
\end{equation*}
Then, from \eqref{sum} and \eqref{eq:caseV}, since the curvature
identity \eqref{main equation} holds for any 4-dimensional manifold,
taking account of \eqref{eq:caseI_relation}, \eqref{eq:re_caseII},
\eqref{eq:re_caseIII}, \eqref{562} and \eqref{10},
we have the following coefficients $c_i$'s :\\
    \begin{equation*}
    \begin{gathered}
    c_1=\frac{\lambda}{4}, \ \ c_2=-\lambda, \ \
    c_3=\frac{\lambda}{4}, \ \ c_4=-\lambda, \ \ c_5=2\lambda,\\
    c_6=\lambda,\ \ \ c_7=-\lambda,\ \ \ c_8=0,\ \ \ c_9=0,\ \ \
    c_{10}=0.
    \end{gathered}
    \end{equation*}

From the above observation, we see that Equation \eqref{main
equation} is unique on a 4-dimensional Riemannian manifold. That is,
the curvature identity \eqref{main equation} for a 4-dimensional
Riemannian manifold is universal.

{\bf{Remark}} The universal relation still holds in the
pseudo-Riemannian setting from the appropriate adjustments of sign
of the metric in the test manifold.
 We refer to \cite{P}.

\end{document}